\documentclass[oneside,english]{amsart}
\usepackage{amsthm}

\makeatletter
\numberwithin{equation}{section} 
\numberwithin{figure}{section} 
\theoremstyle{plain}
\theoremstyle{plain}
\newtheorem{thm}{Theorem}
  \theoremstyle{plain}
  \newtheorem{prop}[thm]{Proposition}
  \theoremstyle{plain}
  \newtheorem{cor}[thm]{Corollary}

\makeatother

\usepackage{babel}

\begin{document}

\title{The Unbounded Commutant of an Operator of Class $C_{0}$}

\author{Hari Bercovici}

\address{Department of Mathematics, Indiana University, Bloomington, IN 47405}

\email{bercovic@indiana.edu}

\thanks{The author was supported in part by a grant from the National Science
Foundation.}
\begin{abstract}
We describe the closed, densely defined linear transformations commuting
with a given operator $T$ of class $C_{0}$ in terms of bounded operators
in $\{T\}'$. Our results extend those of Sarason for operators with
defect index $1$, and Martin in the case of an arbitrary finite defect
index.
\end{abstract}
\maketitle

\section{Introduction}

There has been some interest recently in the study of closed unbounded
linear transformations in the commutant of a bounded operator. For
instance, let $T$ denote the restriction of the backward unilateral
shift to a proper invariant subspace. Then Sarason \cite{sar-oam}
showed that any closed, densely defined linear transformation commuting
with $T$ is of the form $v(T)^{-1}u(T),$ where $u,v\in H^{\infty}$
and $v(T)$ is injective. This extends his earlier result \cite{sar-tams}
pertaining to bounded operators, for which one can take $v=1$.

It is fairly easy to see for the above example that closed linear
transformations commuting with $T$ must in fact commute with every
operator in the commutant $\{T\}'$. Therefore Sarason's theorem can
be viewed as a particular case of a result of Martin \cite{martin},
which we describe next. Assume that $T$ is an operator of class $C_{0}(N)$
as defined in \cite[Chapter III]{SnFbook}, and $X$ is a closed,
densely defined linear transformation commuting with every operator
in $\{T\}'$. Then Martin \cite{martin} proved that $X=v(T)^{-1}u(T)$
with $u,v\in H^{\infty}$ such that $v(T)$ is injective. Thus these
linear transformations are exactly the ones that can be obtained by
applying the Sz.-Nagy---Foias functional calculus \cite[Chapter IV]{SnFbook}
with unbounded functions.

Martin conjectured that his result would be true for operators $T$
of class $C_{0}$ with finite multiplicity. We will show that it is
in fact possible to extend this result to arbitrary contractions of
class $C_{0}$. This follows from a more general description of closed,
densely defined linear transformations $X$ commuting with $T$. In
case $T$ has finite multiplicity, our result states that every such
linear transformation $X$ can be written as $X=v(T)^{-1}Y$, where
$Y$ is a bounded operator in $\{T\}'$, and $v\in H^{\infty}$ is
such that $v(T)$ is injective.

\section{Preliminaries}

We will denote by $\mathcal{B}(\mathcal{H},\mathcal{H}')$ the space
of bounded linear operators $W:\mathcal{H}\to\mathcal{H}'$, where
$\mathcal{H}$ and $\mathcal{H}'$ are complex Hilbert spaces. We
will also write $\mathcal{B}(\mathcal{H})=\mathcal{B}(\mathcal{H},\mathcal{H})$.
Recall that an operator $T\in\mathcal{B}(\mathcal{H})$ is a \emph{quasiaffine
transform} of $T'\in\mathcal{B}(\mathcal{H}')$ if there exists a
\emph{quasiaffinity}, i.e. an injective operator with dense range,
$W\in\mathcal{B}(\mathcal{H},\mathcal{H}')$ satisfying $WT=T'W$.
We write $T\prec T'$ if $T$ is a quasiaffine transform of $T'$.
The operators $T$ and $T'$ are \emph{quasisimilar} if $T\prec T'$
and $T'\prec T$, in which case we write $T\sim T'$.

Assume that $T\in\mathcal{B}(\mathcal{H})$ is a contraction, i.e.
$\|T\|\le1$, and it is completely nonunitary in the sense that it
does not have any nontrivial unitary direct summand. The Sz.-Nagy---Foias
functional calculus \cite[Chapter III]{SnFbook} is an algebra homomorphism
$u\mapsto u(T)\in\mathcal{B}(\mathcal{H})$ of the algebra $H^{\infty}$
of bounded analytic functions in the unit disk, which extends the
usual polynomial calculus. The operator $T$ is said to be of class
$C_{0}$ if $u(T)=0$ for some $u\in H^{\infty}\setminus\{0\}$. When
$T$ is of class $C_{0}$, the ideal $\{u\in H^{\infty}:u(T)=0\}$
is of the form $mH^{\infty}$, where $m$ is an inner function, uniquely
determined up to a constant factor of absolute value $1$, and called
the \emph{minimal function} of $T$. For any inner function $m$,
there exist operators of class $C_{0}$ with minimal function $m$.
The most basic example is constructed as follows. Denote by $S$ the
unilateral shift on the Hardy space $H^{2}$, i.e. $(Sf)(\lambda)=\lambda f(\lambda)$
for $f\in H^{2}$. The space $\mathcal{H}(m)=H^{2}\ominus mH^{2}$
is invariant for $S^{*}$, and the operator $S(m)\in\mathcal{B}(\mathcal{H}(m))$
is defined by the requirement that $S(m)^{*}=S^{*}|\mathcal{H}(m)$.
The operator $S(m)$ has minimal function equal to $m$.

Quasisimilarity allows a complete classification of operators of class
$C_{0}$. We will only need the facts collected in the following statement.
We refer to \cite[Theorem III.5.1]{key-7} for (1-3), \cite[Theorem VII.1.9]{key-7}
for (4), \cite[Proposition III.5.33]{key-7} for (5), \cite[Proposition III.4.7]{SnFbook}
or \cite[Proposition II.4.9]{key-7} for (6), \cite[Proposition VII.1.21]{key-7}
for (7), and \cite[Theorem IV.1.2]{key-7} for (8).
\begin{thm}
\label{thm:grand-preliminary}Let $T\in\mathcal{B}(\mathcal{H})$
and $T'\in\mathcal{B}(\mathcal{H}')$ be operators of class $C_{0}$.
Denote by $m$ the minimal function of $T$.
\begin{enumerate}
\item We have $T\prec T'$ if and only if $T'\prec T$.
\item There exists a collection $\{m_{i}\}_{i\in I}$ of inner divisors
of $m$ such that $m=m_{i}$ for some $i$, and $T\sim\bigoplus_{i\in I}S(m_{i})$.
\item If $T$ has finite cyclic multiplicity $n$, we have $T\sim\bigoplus_{j=1}^{n}S(m_{j})$,
with $m_{1}=m$ and $m_{j+1}$ divides $m_{j}$ for $j=1,2,\dots,n-1$.
\item If $T$ has finite multiplicity, and $\mathcal{M}$ is an invariant
subspace for $T$ such that $T\sim T|\mathcal{M}$, then $\mathcal{M}=\mathcal{H}$.
\item Every invariant subspace $\mathcal{M}$ for $T$ is of the form $\mathcal{M}=\overline{A\mathcal{H}}$,
with $A$ in the commutant $\{T\}'$ of $T$.
\item An operator of the form $v(T)$ with $v\in H^{\infty}$ is injective
if and only if $v$ and $m$ have no nonconstant common inner factors.
In this case, $v(T)$ is a quasiaffinity.
\item If $T$ has finite multiplicity and $A\in\{T\}'$ is injective, then
the map $\mathcal{M}\mapsto\overline{A\mathcal{M}}$ is an inclusion
preserving automorphism of the lattice of invariant subspaces of $T$.
\item For every $Y$ in the double commutant $\{T\}''$ there exist $u,v\in H^{\infty}$
such that $v(T)$ is a quasiaffinity and $Y=v(T)^{-1}u(T)$.
\end{enumerate}
\end{thm}
The following result appears in \cite[Lemma 2.7]{BFsN} (see also
\cite[Proposition IV.1.13]{key-7}), but unfortunately only for multiplicity
2. The argument here follows a different path.
\begin{prop}
\label{pro:AB=00003Dv(T)}Assume that $T\in\mathcal{B}(\mathcal{H})$
is of class $C_{0}$ and has finite multiplicity. For every injective
$A\in\{T\}'$ there exits another injective $B\in\{T\}'$, and a function
$v\in H^{\infty}$ such that $AB=BA=v(T)$. The operators $A,B$ and
$v(T)$ are then quasiaffinities.\end{prop}
\begin{proof}
As seen in \cite{BFsN}, it suffices to consider operators of the
form $T=\bigoplus_{j=1}^{n}S(m_{j}),$ where $m_{j+1}$ divides $m_{j}$
for $j=1,2,\dots,n-1$. Let $A\in\{T\}'$ be an injective operator.
By Theorem \ref{thm:grand-preliminary}(7), the map $\mathcal{M}\mapsto\overline{A\mathcal{M}}$
is an order preserving automorphism of the lattice of invariant subspaces
for $T$. Regard $\mathcal{H}(m_{j})$ as subspaces of $\mathcal{H}=\bigoplus_{j=1}^{n}\mathcal{H}(m_{j})$,
and set $\mathcal{H}_{j}=\overline{A\mathcal{H}(m_{j})}$, $\mathcal{K}_{j}=\bigvee_{i\ne j}\mathcal{H}_{i}$,
and $\mathcal{H}'_{j}=\mathcal{H}\ominus\mathcal{K}_{j}$ for $j=1,2,\dots,n$.
We must then have $\bigcap_{j=1}^{n}\mathcal{K}_{j}=\{0\}$, $\mathcal{H}_{j}\cap\mathcal{K}_{j}=\{0\}$
and $\mathcal{H}_{j}\vee\mathcal{K}_{j}=\mathcal{H}$. The last two
equalities imply that the operator $X_{j}\in\mathcal{B}(\mathcal{H}_{j},\mathcal{H}'_{j})$
defined by $X_{j}=P_{\mathcal{H}'_{j}}|\mathcal{H}_{j}$ is a quasiaffinity.
Moreover, this operator satisfies the equation $X_{j}(T|\mathcal{H}_{j})=T_{j}X_{j}$,
where $T_{j}\in\mathcal{L}(\mathcal{H}'_{j})$ is defined by the equality
$T_{j}^{*}=T^{*}|\mathcal{H}'_{j}$. Thus $T|\mathcal{H}_{j}\prec T_{j}$,
and since $S(m_{j})\prec T|\mathcal{H}_{j}$ (via the operator $A|\mathcal{H}(m_{j})$),
there must exist a quasiaffinity $Y_{j}\in\mathcal{B}(\mathcal{H}'_{j},\mathcal{H}(m_{j}))$
satisfying $Y_{j}T_{j}=S(m_{j})Y_{j}$. We define now an operator
$C\in\{T\}'$ by setting \[
Ch=\bigoplus_{j=1}^{n}Y_{j}P_{\mathcal{H}'_{j}}h.\]
 It is easy to verify that $C$ is a quasiaffinity. Indeed, $Ch=0$
implies that $P_{\mathcal{H}'_{j}}h=0$, and hence $h\in\bigcap_{j=1}^{n}\mathcal{K}_{j}=\{0\}.$
Also, $C\mathcal{H}=\bigvee_{j=1}^{n}Y_{j}\mathcal{H}'_{j}=\mathcal{H}$.
The product $AC$ leaves all the summands $\mathcal{H}(m_{j})$ invariant,
and therefore Sarason's generalized interpolation theorem \cite{sar-tams}
implies the existence of functions $u_{j}\in H^{\infty}$ such that
$AC=\bigoplus_{j=1}^{n}u_{j}(S(m_{j}))$. Moreover, $u_{j}$ and $m_{j}$
have no nonconstant common inner factor because $AC$ is injective.
We deduce from \cite[Theorem III.1.14]{key-7} that there exist scalars
$t_{j}$ such that $v_{j}=u_{j}+t_{j}m_{j}$ has no nonconstant common
inner factor with the minimal function $m_{1}$ of $T$. Note that
we also have $AC=\bigoplus_{j=1}^{n}v_{j}(S(m_{j}))$. Define now
$v=v_{1}v_{2}\cdots v_{n}\in H^{\infty}$ and operators $D,B\in\{T\}'$
by $D=\bigoplus_{j=1}^{n}(v/v_{j})(S(m_{j}))$ and $B=CD$. We have
$AB=v(T)$ and $A(BA-v(T))=ABA-v(T)A=0$ so that $BA=v(T)$ because
$A$ is injective. The operator $v(T)$ is a quasiaffinity because
$v$ and $m_{1}$ do not have nonconstant common inner divisors.
\end{proof}

\section{Unbounded linear transformations in the Commutant}

Consider a Hilbert space $\mathcal{H}$ and a  linear transformation
$X:\mathcal{D}(X)\to\mathcal{H}$, where $\mathcal{D}(X)\subset\mathcal{H}$
is a dense linear manifold. Recall that $X$ is said to be closed
if its graph\[
\mathcal{G}(X)=\{h\oplus Xh:h\in\mathcal{D}(X)\}\]
is a closed subspace in $\mathcal{H}\oplus\mathcal{H}$. The linear
transformation $X$ is closable if the closure $\overline{\mathcal{G}(X)}$
is the graph of a linear transformation, usually denoted $\overline{X}$
and called the closure of $X$.

Let now $T\in\mathcal{B}(\mathcal{H})$ be a completely nonunitary
contraction, let $v\in H^{\infty}$ be such that $v(T)$ is a quasiaffinity,
and let $A\in\{T\}'$. The linear transformation $X=v(T)^{-1}A$ with
domain\[
\mathcal{D}(X)=\{h\in\mathcal{H}:Ah\in v(T)\mathcal{H}\}\]
has graph\[
\mathcal{G}(X)=\{h\oplus k:Ah=v(T)k\},\]
so that $X$ is obviously closed. Moreover, since $v(T)A=Av(T)$,
we have\[
\mathcal{G}(X)\supset\mathcal{G}(Av(T)^{-1})=\{v(T)h\oplus Ah:h\in\mathcal{H}\}\]
and thus $\mathcal{D}(X)\supset v(T)\mathcal{H}$ is dense. If $v_{1}\in H^{\infty}$
is another function such that $v_{1}(T)$ is a quasiaffinity, the
equality $v(T)^{-1}Ah=v_{1}(T)^{-1}A_{1}h$ for $h$ in a dense linear
manifold $\mathcal{D}\subset\mathcal{D}(v(T)^{-1}A)\cap\mathcal{D}(v_{1}(T)^{-1}A_{1})$
implies $v(T)^{-1}A=v_{1}(T)^{-1}A_{1}$. Indeed, we have $v_{1}(T)Ah=v(T)A_{1}h$
for $h\in\mathcal{D}$, hence $v_{1}(T)A=v(T)A_{1}$. Then we deduce
\[
v_{1}(T)((v(T)k-Ah)=v(T)(v_{1}(T)k-A_{1}h),\quad h,k\in\mathcal{H},\]
so that $h\oplus k\in\mathcal{G}(v(T)^{-1}A)$ if and only if $h\oplus k\in\mathcal{G}(v_{1}(T)^{-1}A_{1})$.
These remarks apply more generally to linear transformations of the
form $B^{-1}A$, where $A,B\in\{T\}'$, $B$ is a quasiaffinity, and
$AB=BA$. When $A$ and $B$ do not commute, the linear transformation
$B^{-1}A$ is still closed, but might not be densely defined, while
$AB^{-1}$ is densely defined but perhaps not closable.

Linear transformations of the form $v(T)^{-1}A$, $A\in\{T\}'$, commute
with $T$ in the sense that $TX\subset XT$ or, equivalently, $\mathcal{G}(X)$
is invariant for $T\oplus T$.
\begin{prop}
\label{pro:X=00003DAB-inv}Let $T\in\mathcal{B}(\mathcal{H})$ be
an operator of class $C_{0}$, and let $X$ be a closed, densely defined
linear transformation commuting with $T$. There exist bounded operators
$A,B\in\{T\}'$ such that $B$ is a quasiaffinity and $X=\overline{AB^{-1}}$.\end{prop}
\begin{proof}
The operator $T'=(T\oplus T)|\mathcal{G}(X)$ is of class $C_{0}$,
and $T'\prec T$. Indeed, the operator $W\in\mathcal{B}(\mathcal{G}(X),\mathcal{H})$
defined by $W(h\oplus k)=h$ satisfies $WT'=TW$, and $W$ is injective
(because $\mathcal{G}(X)$ is a graph) and has dense range $\mathcal{D}(X)$.
Theorem \ref{thm:grand-preliminary}(1) implies the existence of an
injective operator $V\in\mathcal{B}(\mathcal{H},\mathcal{H}\oplus\mathcal{H})$
such that $\overline{V\mathcal{H}}=\mathcal{G}(X)$ and $(T\oplus T)V=VT$.
Writing $Vh=Bh\oplus Ah$ for $h\in\mathcal{H}$, the operators $A,B$
must belong to $\{T\}'$. Moreover, $B$ is a quasiaffinity. Indeed,
$Bh=0$ implies $Ah=XBh=0$, so that $Vh=0$ and hence $h=0$ because
$V$ is injective. The fact that $V\mathcal{H}$ is dense in $\mathcal{G}(X)$
implies that $\overline{B\mathcal{H}}\supset\mathcal{D}(X)$, and
hence $B$ has dense range. Obviously $\mathcal{G}(AB^{-1})=V\mathcal{H}$,
and hence $X=\overline{AB^{-1}}$.
\end{proof}
For operators with finite multiplicity, a stronger result can be proved.
\begin{thm}
\label{thm:finite-multiplicity} Let $T\in\mathcal{B}(\mathcal{H})$
be an operator of class $C_{0}$ with finite multiplicity, and let
$X$ be a closed, densely defined linear transformation commuting
with $T$. There exist $A\in\{T\}'$ and $v\in H^{\infty}$ such that
$v(T)$ is a quasiaffinity and $X=v(T)^{-1}A$.\end{thm}
\begin{proof}
By Proposition \ref{pro:X=00003DAB-inv}, we can find $A_{0},B\in\{T\}'$
such that $B$ is a quasiaffinity and $X\supset A_{0}B^{-1}$. Proposition
\ref{pro:AB=00003Dv(T)} implies the existence of $v\in H^{\infty}$
and of a quasiaffinity $C\in\{T\}'$ such that $BC=CB=v(T).$ Setting
now $A=A_{0}C$, we have \[
Av(T)^{-1}=A_{0}C(BC)^{-1}\subset A_{0}B^{-1}\subset X.\]
We conclude the proof by showing that both $v(T)^{-1}A$ and $X$
coincide with the closure of $Av(T)^{-1}$. For this purpose, define
operators $T_{1}=(T\oplus T)|\mathcal{G}(X)$, $T_{2}=(T\oplus T)|\mathcal{G}(v(T)^{-1}A)$,
and $T_{3}=(T\oplus T)|\mathcal{G}(\overline{Av(T)^{-1}})$. As observed
earlier, $T_{1}\sim T_{2}\sim T_{3}\sim T$. Since $\mathcal{G}(\overline{Av(T)^{-1}})$
is an invariant subspace for $T_{1}$ and $T_{2}$, theorem \ref{thm:grand-preliminary}(4)
implies the desired conclusion that $X=v(T)^{-1}A$.
\end{proof}
Our final result pertains to double commutants.
\begin{thm}
\label{thm:double-comm}Let $T\in\mathcal{B}(\mathcal{H})$ be an
operator of class $C_{0}$, and let $X$ be a closed, densely defined
linear transformation commuting with every $A\in\{T\}'$. Then there
exist $u,v\in H^{\infty}$ such that $v(T)$ is a quasiaffinity and
$X=v(T)^{-1}u(T)$.\end{thm}
\begin{proof}
We first prove the result under the additional assumption that $T$
has finite multiplicity. In this case, Theorem \ref{thm:finite-multiplicity}
yields $A_{0}\in\{T\}'$ and $v_{0}\in H^{\infty}$ such that $v_{0}(T)$
is a quasiaffinity and $X=v_{0}(T)^{-1}A_{0}$. We observe next that
$A_{0}$ belongs to the double commutant $\{T\}''$. Indeed, for any
$B\in\{T\}'$ and $h\in\mathcal{D}(X)$ we have $Bh\in\mathcal{D}(X)$
and $XBh=BXh$ so that\[
v_{0}(T)XBH=v_{0}(T)BXH=Bv_{0}(T)Xh\]
 and therefore $A_{0}Bh=BA_{0}h$. We conclude that $A_{0}B=BA_{0}$
because $\mathcal{D}(X)$ is dense. By Theorem \ref{thm:grand-preliminary}(8),
there exist $u,v_{1}\in H^{\infty}$ such that $v_{1}(T)$ is a quasiaffinity
and $A_{0}=v_{1}(T)^{-1}u(T).$ We reach the desired conclusion $X=v(T)^{-1}u(T)$
with $v=v_{0}v_{1}$.

Consider now an arbitrary operator of class $C_{0}$, and let $m$
denote its minimal function. Let $\mathcal{M}\subset\mathcal{H}$
be an invariant subspace for $T$ such that $T|\mathcal{M}$ has finite
multiplicity and minimal function equal to $m$. By Theorem \ref{thm:grand-preliminary}(5),
$\mathcal{M}=\overline{C\mathcal{H}}$ for some $C\in\{T\}'$. We
have $C\mathcal{D}(X)\subset\mathcal{D}(X)\cap\mathcal{M}$ and \[
X(C\mathcal{D}(X))\subset CX\mathcal{D}(X)\subset C\mathcal{H}\subset\mathcal{M}.\]
Therefore there exists a closed densely defined linear transformation
$X_{\mathcal{M}}$ on $\mathcal{M}$ such that\[
\mathcal{G}(X_{\mathcal{M}})=\mathcal{G}(X)\cap(\mathcal{M}\oplus\mathcal{M}).\]
We claim that $\mathcal{D}(X_{\mathcal{M}})=\mathcal{D}(X)\cap\mathcal{M}$.
Indeed, let us set $T_{1}=(T\oplus T)|\mathcal{G}(X_{\mathcal{M}})$
and $T_{2}=(T\oplus T)|\mathcal{G}(X)\cap(\mathcal{M}\oplus\mathcal{H})$.
The projection on the first component demonstrates the relations $T_{1}\prec T|\mathcal{M}$
and $T_{2}\prec T|\mathcal{M}$. The equality\[
\mathcal{G}(X_{\mathcal{M}})=\mathcal{G}(X)\cap(\mathcal{M}\oplus\mathcal{H}),\]
and hence $\mathcal{D}(X_{\mathcal{M}})=\mathcal{D}(X)\cap\mathcal{M}$,
follows from Theorem \ref{thm:grand-preliminary}(4). A similar argument
shows that $\mathcal{G}(X_{\mathcal{M}})$ is the closure of $\{Ch\oplus CXh:h\in\mathcal{D}(X)\}$.

We show next that $X_{\mathcal{M}}$ commutes with every operator
in the commutant of $T|\mathcal{M}$. Indeed, let $D\in\mathcal{B}(\mathcal{M})$
be such an operator. Then $DC\in\{T\}'$ so that $DCh\in\mathcal{D}(X)$
for every $h\in\mathcal{D}(X)$, and\[
XDCh=DCXh=DXCh.\]
Thus $D\oplus D$ leaves $\{Ch\oplus CXh:h\in\mathcal{D}(X)\}$ invariant,
and hence it leaves its closure invariant as well, i.e. $D$ commutes
with $X_{\mathcal{M}}$.

The first part of the proof implies the existence of $u,v\in H^{\infty}$
such that $v(T_{\mathcal{M}})$ is a quasiaffinity, and $X_{\mathcal{M}}=v(T|\mathcal{M})^{-1}u(T|\mathcal{M})$.
Note that $v(T)$ is a quasiaffinity as well since $T$ and $T|\mathcal{M}$
have the same minimal function (cf. Theorem \ref{thm:grand-preliminary}(6)).
We claim that $X=v(T)^{-1}u(T)$. Indeed, consider arbitrary vectors
$h_{1}\in\mathcal{D}(X)$, $h_{2}\in\mathcal{D}(v(T)^{-1}u(T))$,
and let $\mathcal{M}_{1}\supset\mathcal{M}$ be an invariant subspace
for $T$ such that $T|\mathcal{M}_{1}$ has finite multiplicity, and
$h_{1},h_{2}\in\mathcal{M}_{1}$; for instance, once can take $\mathcal{M}_{1}$
to be the smallest invariant subspace containing $\mathcal{M},h_{1}$
and $h_{2}$. The preceding argument, with $\mathcal{M}_{1}$ in place
of $\mathcal{M}$, shows that $X_{\mathcal{M}_{1}}=v_{1}(T|\mathcal{M}_{1})^{-1}u_{1}(T|\mathcal{M}_{1})$
for some $u_{1},v_{1}\in H^{\infty}$ such that $v_{1}(T)$ is a quasiaffinity.
Note now that, for $h\in\mathcal{D}(X)\cap\mathcal{M}$, we have both
$v(T)Xh=u(T)h$ and $v_{1}(T)Xh=u_{1}(T)h,$ and therefore\[
(v_{1}(T)u(T)-v(T)u_{1}(T))h=v_{1}(T)v(T)Xh-v(T)v_{1}(T)Xh=0\]
for such vectors. Since $\mathcal{D}(X)\cap\mathcal{M}$ is dense
in $\mathcal{M}$, we have $(v_{1}u-u_{1}v)(T|\mathcal{M})=0$. We
deduce that $m$, which is the minimal function of $T|\mathcal{M}$,
divides $v_{1}u-vu_{1}$, and thus $v_{1}(T)u(T)=v(T)u_{1}(T)$. This
implies that $v(T)^{-1}u(T)=v_{1}(T)^{-1}u_{1}(T)$, and therefore\[
h_{1}\in\mathcal{D}(X)\cap\mathcal{M}_{1}=\mathcal{D}(X_{\mathcal{M}_{1}})=\mathcal{D}(v_{1}(T|\mathcal{M}_{1})^{-1}u(T|\mathcal{M}_{1}))\subset\mathcal{D}(v(T)^{-1}u(T)),\]
\begin{eqnarray*}
h_{2} & \in & \mathcal{D}(v(T)^{-1}u(T))\cap\mathcal{M}_{1}=\mathcal{D}(v(T|\mathcal{M}_{1})^{-1}u(T|\mathcal{M}_{1}))\\
 & = & \mathcal{D}(v_{1}(T|\mathcal{M}_{1})^{-1}u(T|\mathcal{M}_{1}))=\mathcal{D}(X_{\mathcal{M}_{1}})\subset\mathcal{D}(X),\end{eqnarray*}
and\[
Xh_{j}=v_{1}(T)^{-1}u_{1}(T)h_{j}=v(T)^{-1}u(T)h_{j}\]
for $j=1,2$. The desired equality $X=v(T)^{-1}u(T)$ follows.
\end{proof}
When $T$ has multiplicity $1$, i.e. $T$ has a cyclic vector, the
algebra $\{T\}'$ is precisely the algebra genreated by $T$ and closed
in the weak operator topology; see \cite[Theorem IV.1.2]{key-7}.
Therefore Theorem \ref{thm:double-comm} implies the following extension
of Sarason's result \cite{sar-oam}.
\begin{cor}
Let $T\in\mathcal{B}(\mathcal{H})$ be an operator of class $C_{0}$
with multiplicity $1$, and let $X$ be a closed, densely defined
linear transformation commuting with $T$. Then there exist $u,v\in H^{\infty}$
such that $v(T)$ is a quasiaffinity and $X=v(T)^{-1}u(T)$.
\end{cor}
In Theorem \ref{thm:double-comm}, if  we only assume that $X$ is
a densely defined linear transformation commuting with $\{T\}'$,
the conclusion is that $X\subset v(T)^{-1}u(T)$ for some $u,v\in H^{\infty}$
such that $v(T)$ is a quasiaffinity. Indeed, the operator $X$ must
be closable by \cite[Proposition 5.8]{closability}. As noted by Martin,
in case $T=S(m)$ this was also proved by Sarason \cite[Lemma 3]{martin}.

\end{document}